\documentclass[reqno,10pt]{amsart}

\usepackage{amssymb}
\usepackage[applemac]{inputenc}
\usepackage[T1]{fontenc}

\usepackage[english,frenchb]{babel}

\newtheorem{theorem}{Theorem}[section]

\newtheorem{e-proposition}[theorem]{Proposition}

\newtheorem{e-definition}[theorem]{Definition\rm}

\newtheorem{theoreme}{Th\'eor\`eme}[section]

\begin{document}

\selectlanguage{english}
\title[Perturbation Method for the Einstein-Dirac-Maxwell Equations]{Perturbation Method for Particle-like Solutions of the Einstein-Dirac-Maxwell Equations}

\selectlanguage{english}
\author{Simona Rota Nodari}
\email{rotanodari@ceremade.dauphine.fr}

\address{Ceremade (UMR CNRS 7534),
Universit\'e Paris-Dauphine,
Place du Mar\'echal de Lattre de Tassigny,
75775 Paris Cedex 16,
France}
\address{Dipartimento di Matematica, Universit\`a di Milano, Via Saldini 50, 20133 Milan, Italy}

\date{December 18, 2009}

\begin{abstract}
\selectlanguage{english}

The aim of this Note is to prove by a perturbation method the existence of solutions of the coupled Einstein-Dirac-Maxwell equations for a static, spherically symmetric system of two fermions in a singlet spinor state and with the electromagnetic coupling constant $\left(\frac{e}{m}\right)^2<1$. 
We show that the nondegenerate solution of Choquard's equation generates a branch of solutions of the Einstein-Dirac-Maxwell equations.

\medskip

\selectlanguage{frenchb}
 
\noindent {\bf Une m\'ethode de perturbation pour les solutions localis\'ees des \'equa\-tions d'Einstein-Dirac-Maxwell.}

\medskip

\noindent\textsc{R\'esum\'e.} 
Le but de cette Note est de d\'emontrer par une m\'ethode de perturbation l'existence de solutions des \'equations d'Einstein-Dirac-Maxwell pour un syst\`eme statique, \`a sym\'etrie sph\'erique  de deux fermions dans un \'etat de singulet et avec une constante de couplage \'electromagn\'etique $\left(\frac{e}{m}\right)^2<1$. 
On montre que la solution non d\'eg\'en\'er\'ee de l'\'equation de Choquard g\'en\`ere une branche de solutions des \'equations d'Einstein-Dirac-Maxwell. 

\end{abstract}

\maketitle

\selectlanguage{frenchb}
\section*{Version fran\c{c}aise abr\'eg\'ee}

Dans un papier r\'ecent \cite{rotanodari}, par une m\'ethode de perturbation, on a montr\'e de mani\`ere rigoureuse l'existence de solutions des \'equations d'Einstein-Dirac pour un syst\`eme statique, \`a sym\'etrie sph\'erique  de deux fermions dans un \'etat de singulet. Dans cette Note, on g\'en\'eralise ce r\'esultat aux \'equations d'Einstein-Dirac-Maxwell et on montre, dans le cas particulier d'un couplage \'electromagn\'etique faible, l'existence des solutions obtenues num\'eriquement par F. Finster, J. Smoller et ST. Yau dans \cite{finsmoyaumax}.

Plus pr\'ecis\'ement, en utilisant l'id\'ee introduite par Ounaies pour une classe d'\'equations de Dirac non lin\'eaires (voir  \cite{ounaies}) et adapt\'ee dans  \cite{rotanodari} aux \'equations d'Einstein-Dirac, on obtient le th\'eor\`eme suivant. 

\begin{theoreme}\label{th:solutionfr} Soient $e,m,\omega$ tels que $e^2-m^2<0$, $0<\omega< m$ et supposons $m-\omega$ assez petit; alors il existe une solution non triviale de (\ref{eq:einsteindiracmaxwelleq1}-\ref{eq:einsteindiracmaxwelleq5}).
\end{theoreme}

Dans cette Note, on d\'ecrit la m\'ethode utilis\'ee pour d\'emontrer ce th\'eor\`eme. Premi\`erement, par un changement d'\'echelle, on transforme les \'equations d'Einstein-Dirac-Maxwell (\ref{eq:einsteindiracmaxwelleq2.1}-\ref{eq:einsteindiracmaxwelleq2.4}) en un syst\`eme perturb\'e qui s'\'ecrit sous la forme (\ref{eq:systemepsilon}). On choisit $\varepsilon=m-\omega$ comme param\`etre de perturbation.

Deuxi\`ement, on remarque que, pour $\varepsilon=0$ et $\left(\frac{e}{m}\right)^2<1$, ce syst\`eme est \'equivalent au syst\`eme (\ref{eq:systemepsilonzero2}) o\`u l'\'equation pour la variable $\varphi$ est l'\'equation de Choquard. Il est bien connu que l'\'equation de Choquard a une solution radiale positive. De plus, dans l'espace des fonctions radiales, cette solution est non d\'eg\'en\'er\'ee, dans le sens o\`u le noyau de la lin\'earisation de l'\'equation contient seulement la fonction identiquement nulle.
On appelle $\phi_0$ la solution du syst\`eme (\ref{eq:systemepsilonzero2}).

Ensuite, on observe que le syst\`eme perturb\'e s'\'ecrit sous la forme $D(\varepsilon,\varphi,\chi,\tau,\zeta)=0$ avec  $D$ un op\'erateur non lin\'eaire de classe $\mathcal{C}^1$, pour un bon choix d'espaces fonctionnels. On prouve que cet op\'erateur satisfait les hypoth\`eses du th\'eor\`eme des fonctions implicites. En particulier, on montre que la lin\'earisation de l'op\'erateur $D$ par rapport \`a $(\varphi,\chi,\tau,\zeta)$ en $(0,\phi_0)$, $D_{\varphi,\chi,\tau,\zeta}(0,\phi_0)$, est une injection, gr\^ace \`a la non-d\'eg\'en\'erescence de la solution de l'\'equation de Choquard, et s'\'ecrit comme somme d'un isomorphisme et d'un op\'erateur compact; donc $D_{\varphi,\chi,\tau,\zeta}(0,\phi_0)$  est un isomorphisme. En appliquant  le  th\'eor\`eme des fonctions implicites, on d\'eduit que, pour $\varepsilon$ assez petit et $e^2-m^2<0$, le syst\`eme (\ref{eq:systemepsilon}) a une solution.

 En conclusion, pour $e^2-m^2<0$, $0<\omega< m$ et $m-\omega$ assez petit, les \'equations d'Einstein-Dirac-Maxwell poss\`edent  une solution non triviale.

\selectlanguage{english}

\section{Introduction}
In a recent paper \cite{rotanodari}, using a perturbation method, we proved rigorously the existence of solutions of the coupled Einstein-Dirac equations for a static, spherically symmetric system of two fermions in a singlet spinor state. In this Note, we extend our result to the Einstein-Dirac-Maxwell equations and we prove, in the particular case of a weak electromagnetic coupling, the existence of the solutions obtained numerically by F. Finster, J. Smoller and ST. Yau in \cite{finsmoyaumax}. 

The general Einstein-Dirac-Maxwell equations for a system of $n$ Dirac particles take the form
$$
\begin{array}{lll}
(G-m)\psi_a=0, & R^i_j-\frac{1}{2}R\delta^i_j=-8\pi T^i_j, &\nabla_kF^{jk}=4\pi e\sum\limits_{a=1}^n\overline{\psi_a}G^j\psi_a
\end{array}
$$
where $G^j$ are the Dirac matrices, $G$ denote the Dirac operator, $\psi_a$ are the wave functions of fermions of mass $m$ and charge $e$, $F_{jk}$ is the electromagnetic field tensor and, finally, $T^i_j$ is the sum of the energy-momentum tensor of the Dirac particle and the Maxwell stress-energy tensor.

In \cite{finsmoyaumax}, the metric, in polar coordinates $(t,r,\vartheta,\varphi)$, is given by $$ds^2=\frac{1}{T^2}\,dt^2-\frac{1}{A}\,dr^2-r^2\,d\vartheta^2-r^2\sin^2\vartheta\,d\varphi^2$$
with $A=A(r)$, $T=T(r)$ positive functions; moreover, using the ansatz from \cite{finsmoyau}, Finster, Smoller and Yau describe the Dirac spinors with two real radial functions 
$\Phi_1(r)$, $\Phi_2(r)$
and they assume that the electromagnetic potential has the form $\mathcal{A}=(-V,0)$, with  $V$ the Coulomb potential.

In this case the Einstein-Dirac-Maxwell equations can be written as 
\begin{eqnarray}\label{eq:einsteindiracmaxwelleq1}
\sqrt{A}\Phi_1'&=&\frac{1}{r}\Phi_1-((\omega-eV) T+m)\Phi_2\\
\label{eq:einsteindiracmaxwelleq2}
\sqrt{A}\Phi_2'&=&((\omega-eV)T-m)\Phi_1-\frac{1}{r}\Phi_2\\
\label{eq:einsteindiracmaxwelleq3}
rA'&=&1-A-16\pi(\omega-eV) T^2\left(\Phi_1^2+\Phi_2^2\right)-r^2AT^2\left(V'\right)^2\\
\label{eq:einsteindiracmaxwelleq4}
2rA\frac{T'}{T}&=&A-1-16\pi(\omega-eV) T^2\left(\Phi_1^2+\Phi_2^2\right)+32\pi\frac{1}{r} T\Phi_1\Phi_2\nonumber\\
&&+16\pi m T\left(\Phi_1^2-\Phi_2^2\right)+r^2AT^2\left(V'\right)^2\\
\label{eq:einsteindiracmaxwelleq5}
r^2AV''&=&-8\pi e \left(\Phi_1^2+\Phi_2^2\right)- \left(2rA+r^2A\frac{T'}{T}+\frac{r^2}{2}A'\right)V'
\end{eqnarray}
with the normalization condition $\int_0^{\infty}{|\Phi|^2\frac{T}{\sqrt A}\,dr}=\frac{1}{4\pi}$. \\
In order that the metric be asymptotically Minkowskian and the solutions have finite (ADM) mass, Finster, Smoller and Yau assume 
$$\lim\limits_{r\rightarrow \infty}T(r)=1$$ and $$\lim\limits_{r\rightarrow \infty}\frac{r}{2}(1-A(r))<\infty.$$
Finally, they also require that the electromagnetic potential vanishes at infinity.

In this Note, using the idea introduced by Ounaies for a class of nonlinear Dirac equations (see \cite{ounaies}) and adapted in \cite{rotanodari} to the Einstein-Dirac equations, we obtain the following result.
\begin{theorem}\label{th:solution} Given $e,m,\omega$ such that $e^2-m^2<0$, $0<\omega< m$ and $m-\omega$ is sufficiently small, there exists a non trivial solution of (\ref{eq:einsteindiracmaxwelleq1}-\ref{eq:einsteindiracmaxwelleq5}).
\end{theorem}

\section{Perturbation method for the Einstein-Dirac-Maxwell equations}\label{section:EinsteinDiracMaxwell}

First of all, we observe that, writing $T(r)=1+t(r)$ and integrating the equation (\ref{eq:einsteindiracmaxwelleq5}), the Einstein-Dirac-Maxwell equations become
\begin{eqnarray}\label{eq:einsteindiracmaxwelleq2.1}
\sqrt{A}\Phi_1'&=&\frac{1}{r}\Phi_1-((\omega-eV)(1+t)+m)\Phi_2\\
\label{eq:einsteindiracmaxwelleq2.2}
\sqrt{A}\Phi_2'&=&((\omega-eV)(1+t)-m)\Phi_1-\frac{1}{r}\Phi_2\\
\label{eq:einsteindiracmaxwelleq2.3}
2rAt'&=&(A-1)(1+t)-16\pi(\omega-eV) (1+t)^3\left(\Phi_1^2+\Phi_2^2\right)\nonumber\\
&&+32\pi\frac{1}{r} (1+t)^2\Phi_1\Phi_2+16\pi m (1+t)^2\left(\Phi_1^2-\Phi_2^2\right)\nonumber\\
&&+r^2A(1+t)^3\left(V'\right)^2\\
\label{eq:einsteindiracmaxwelleq2.4}
\sqrt A (1+t) V'&=&-\frac{8\pi e}{r^2} \int_0^r \left(\Phi_1^2+\Phi_2^2\right) \frac{(1+t)}{\sqrt A}\,ds.
\end{eqnarray}
where $A(r)=1+a(r)$ and
\begin{eqnarray}\label{eq:defa}
a(r)&=&-\frac{1}{r}\exp{(-F(r))}\int_0^{r}\left[16\pi(\omega-eV) (1+t)^2\left(\Phi_1^2+\Phi_2^2\right)\right.\nonumber\\
&&\left.+s^2(1+t)^2\left(V'\right)^2\right]\exp{(F(s))}\,ds
\end{eqnarray}
with $F(r)=\int_0^r s(1+t)^2\left(V'\right)^2\,ds$.

After that, we introduce the new variable $(\varphi,\chi,\tau,\zeta)$ such that 
$$
\begin{array}{llll}
\Phi_1(r)=\varepsilon^{1/2}\varphi(\varepsilon^{1/2} r), &\Phi_2(r)=\varepsilon\chi(\varepsilon^{1/2} r), &t(r)=\varepsilon\tau(\varepsilon^{1/2}r), &V(r)=\varepsilon \zeta(\varepsilon^{1/2}r)
\end{array}
$$
 where $\Phi_1,\Phi_2, t,V$ satisfy  (\ref{eq:einsteindiracmaxwelleq2.1}-\ref{eq:einsteindiracmaxwelleq2.4}) and $\varepsilon=m-\omega$. Using the explicit expression of $a(r)$, given in (\ref{eq:defa}), we write 
 $$a(\Phi_1,\Phi_2,t,V)=\varepsilon \alpha(\varepsilon,\varphi,\chi,\tau,\zeta)$$
 with $\alpha(0,\varphi,\chi,\tau,\zeta)=-\frac{16\pi m}{r}\int_0^{r} \varphi^2\,ds$ . It is now clear that if $\Phi_1,\Phi_2,t,V$ satisfy (\ref{eq:einsteindiracmaxwelleq2.1}-\ref{eq:einsteindiracmaxwelleq2.4}), then $\varphi,\chi,\tau,\zeta$ satisfy the system 
\begin{equation}\label{eq:systemepsilon}
\left\{
\begin{array}{l}
(1+\varepsilon\alpha(\varepsilon,\varphi,\chi,\tau,\zeta))^{1/2}\frac{d}{dr}\varphi-\frac{1}{r}\varphi+2m\chi+K_1\left(\varepsilon,\varphi,\chi,\tau,\zeta \right)=0\\[5pt]
(1+\varepsilon\alpha(\varepsilon,\varphi,\chi,\tau,\zeta))^{1/2} \frac{d}{dr}\chi+\frac{1}{r}\chi+\varphi-m\varphi\tau+e\varphi\zeta+K_2\left(\varepsilon,\varphi,\chi,\tau,\zeta \right)=0\\[5pt]
(1+\varepsilon\alpha(\varepsilon,\varphi,\chi,\tau,\zeta))\frac{d}{dr}\tau-\frac{\alpha(\varepsilon,\varphi,\chi,\tau,\zeta)}{2r}+K_3\left(\varepsilon,\varphi,\chi,\tau,\zeta \right)=0\\[5pt]
(1+\varepsilon\alpha(\varepsilon,\varphi,\chi,\tau,\zeta))^{1/2}(1+\varepsilon\tau)\frac{d}{dr}\zeta+\frac{8\pi e}{r^2}\int_0^r\varphi^2\,ds+K_4\left(\varepsilon,\varphi,\chi,\tau,\zeta \right)=0
\end{array}
\right.
\end{equation}
where $K_1\left(0,\varphi,\chi,\tau,\zeta \right) = K_2\left(0,\varphi,\chi,\tau,\zeta \right) =K_3\left(0,\varphi,\chi,\tau,\zeta \right) =K_4\left(0,\varphi,\chi,\tau,\zeta \right)=0 $. 

Then, for $\varepsilon=0$, (\ref{eq:systemepsilon}) becomes
\begin{equation}\label{eq:systemepsilonzero2}
\left\{
\begin{array}{l}
-\frac{d^2}{dr^2}\varphi+2m\varphi+16\pi (e^2-m^2)m\left(\int_{0}^{\infty}\frac{\varphi^2}{\max(r,s)}\,ds\right)\varphi=0\\[5pt]
\chi(r)=\frac{1}{2m}\left(\frac{1}{r}\varphi-\frac{d}{dr}\varphi\right)\\[5pt] \tau(r)=8\pi m \int_{0}^{\infty}\frac{\varphi^2}{\max(r,s)}\,ds\\[5pt] \zeta(r)=8\pi e \int_{0}^{\infty}\frac{\varphi^2}{\max(r,s)}\,ds\\[5pt]
\end{array}
\right.
\end{equation}
We remark that if $e^2-m^2<0$ the first equation of the system (\ref{eq:systemepsilonzero2}) is the Choquard equation 
\begin{equation}\label{eq:choquard}
\begin{array}{cc}
-\triangle u+2m u-4(m^2-e^2)m\left(\int_{\mathbb{R}^3}\frac{\left|u(y)\right|^2}{|x-y|}\,dy\right)u=0&\mbox{in}\ H^1\left(\mathbb{R}^3\right)
\end{array}
\end{equation}
with $u(x)=\frac{\varphi(|x|)}{|x|}$.
It is well known that Choquard's equation (\ref{eq:choquard}) has a unique radial, positive solution $u_0$ with $\int|u_0|^2=N$ for some $N>0$ given. Furthermore, $u_0$ is infinitely differentiable, goes to zero at infinity and is a radial nondegenerate solution; by this we mean that the linearization of (\ref{eq:choquard}) around $u_0$ has a trivial nullspace in $L^2_{r}(\mathbb{R}^3)$
(see \cite{liebcho}, \cite{Lionscho}, \cite{Lenzmann} for more details). \\
Let $\phi_0=(\varphi_0,\chi_0,\tau_0,\zeta_0)$ be the ground state solution of (\ref{eq:systemepsilonzero2}). 

The main idea is that the solutions of (\ref{eq:systemepsilon}) are the zeros of a $\mathcal{C}^1$ operator $D:\mathbb{R}\times X_{\varphi}\times X_{\chi}\times X_{\tau}\times X_{\zeta}\rightarrow Y_{\varphi}\times Y_{\chi}\times Y_{\tau}\times Y_{\zeta}$. 
So, to obtain a solution of  (\ref{eq:systemepsilon}) from $\phi_0$, we define the operators  
\begin{eqnarray*}
L_1(\varepsilon,\varphi,\chi,\tau,\zeta)&=&(1+\varepsilon\alpha(\varepsilon,\varphi,\chi,\tau,\zeta))^{1/2}\frac{1}{r}\frac{d}{dr}\varphi-\frac{\varphi}{r^2}+2m\frac{\chi}{r}\\&&+\frac{1}{r}K_1\left(\varepsilon,\varphi,\chi,\tau,\zeta \right)\\
L_2(\varepsilon,\varphi,\chi,\tau,\zeta)&=&(1+\varepsilon\alpha(\varepsilon,\varphi,\chi,\tau,\zeta))^{1/2}\frac{1}{r} \frac{d}{dr}\chi+\frac{\chi}{r^2}+\frac{\varphi}{r}-m\frac{\varphi}{r}\tau+e\frac{\varphi}{r}\zeta\\&&+\frac{1}{r}K_2\left(\varepsilon,\varphi,\chi,\tau,\zeta \right)\\
L_3(\varepsilon,\varphi,\chi,\tau,\zeta)&=&(1+\varepsilon\alpha(\varepsilon,\varphi,\chi,\tau,\zeta))\frac{d}{dr}\tau-\frac{\alpha(\varepsilon,\varphi,\chi,\tau,\zeta)}{2r}+K_3\left(\varepsilon,\varphi,\chi,\tau,\zeta \right)\\
L_4(\varepsilon,\varphi,\chi,\tau,\zeta)&=&(1+\varepsilon\alpha(\varepsilon,\varphi,\chi,\tau,\zeta))^{1/2}(1+\varepsilon\tau)\frac{d}{dr}\zeta+\frac{8\pi e}{r^2}\int_0^r\varphi^2\,ds\\&&+K_4\left(\varepsilon,\varphi,\chi,\tau,\zeta \right)
\end{eqnarray*}
and 
$$
D(\varepsilon,\varphi,\chi,\tau,\zeta)=\left(L_1(\varepsilon,\varphi,\chi,\tau,\zeta), L_2(\varepsilon,\varphi,\chi,\tau,\zeta),L_3(\varepsilon,\varphi,\chi,\tau,\zeta),L_4(\varepsilon,\varphi,\chi,\tau,\zeta)\right),
$$
with $X_{\varphi}$, $X_{\chi}$, $X_{\tau}$, $Y_{\varphi}, Y_{\chi}, Y_{\tau}$ defined as in \cite{rotanodari} and 
\begin{eqnarray*}
X_\zeta&=&\left\{\zeta:(0,\infty)\rightarrow\mathbb{R}\left|\lim_{r\rightarrow\infty}\zeta(r)=0, \frac{d}{dr}\zeta\in L^1((0,\infty),dr)\cap L^2((0,\infty),rdr )\right.\right\}\\
Y_\zeta&=&L^1((0,\infty),dr)\cap L^2((0,\infty),rdr )
\end{eqnarray*}
with their natural norms.

Next, we linearize the operator $D$ on $(\varphi,\chi,\tau,\zeta)$ around $(0,\phi_0)$:
\begin{eqnarray*}
D_{\varphi,\chi,\tau,\zeta}(0,\phi_0)(h,k,l,z)&=&\left(
\begin{array}{c}
\frac{1}{r}\frac{d}{dr}h-\frac{h}{r^2}+2m\frac{k}{r}\\[5pt]
\frac{1}{r}\frac{d}{dr}k+\frac{k}{r^2}+\frac{h}{r}-m\frac{\varphi_0}{r} l+e\frac{\varphi_0}{r} z\\[5pt]
\frac{d}{dr}l\\[5pt]
\frac{d}{dr}z
\end{array}
\right)\\
&&+
\left(
\begin{array}{c}
0\\[5pt]
-m \frac{h}{r}\tau_0+e \frac{h}{r}\zeta_0\\[5pt]
\frac{16\pi m}{r^2}\int_0^r\varphi_0h\,ds\\[5pt]
\frac{16\pi e}{r^2}\int_0^r\varphi_0h\,ds\\
\end{array}
\right).
\end{eqnarray*}
We observe that, thanks to the nondegeneracy of the solution of Choquard's equation, $D_{\varphi,\chi,\tau,\zeta}(0,\phi_0)$ is a one-to-one operator. Moreover, it can be written as a sum of an isomorphism and a compact operator. It is thus an isomorphism. Finally, the application of the implicit function theorem yields the following result, which is equivalent to Theorem \ref{th:solution}.
\begin{theorem}\label{th:principalth1} Suppose $e^2-m^2<0$ and let $\phi_0$ be the ground state solution of (\ref{eq:systemepsilonzero2}), then there exists $\delta>0$ and a function $\eta\in\mathcal{C}((0,\delta),X_\varphi\times 
X_\chi\times X_\tau\times X_\zeta)$ such that $\eta(0)=\phi_0$ and $D(\varepsilon,\eta(\varepsilon))=0$ for $0\leq \varepsilon <\delta$.
\end{theorem}

\section*{Acknowledgements} The author would like to thank professor Eric S\'er\'e and professor Bernhard Ruf for helpful discussions.

\end{document}